\documentclass[journal,transmag]{IEEEtran}
\usepackage{amsmath}
\usepackage{hyperref,color,breqn,multirow}
\usepackage[top=0.7in, left=0.65in]{geometry}
\usepackage{graphicx}
\usepackage[printonlyused,withpage,nolist]{acronym}
\usepackage{scalerel}
\usepackage{soul}

\newcommand{\opdiv}{\operatorname{div}}
\newcommand{\opcurl}{\operatorname{curl}}
\newcommand{\veczero}{\mathbf{0}}
\newcommand{\vecB}{\mathbf{B}}
\newcommand{\vecE}{\mathbf{E}}
\newcommand{\vecH}{\mathbf{H}}
\newcommand{\vecJ}{\mathbf{J}}
\newcommand{\vecT}{\mathbf{T}}
\newcommand{\vecV}{\mathbf{V}}
\newcommand{\vecv}{\mathbf{v}}

\definecolor{mygray}{gray}{.5}

\begin{document}
\title{An Equilibrated Error Estimator for the 2D/1D MSFEM T-Formulation of the Eddy Current Problem}
\author{\IEEEauthorblockN{Markus Sch\"obinger\IEEEauthorrefmark{1},
Karl Hollaus\IEEEauthorrefmark{1}
}
\IEEEauthorblockA{\IEEEauthorrefmark{1}Technische Universit\"at Wien, Institute for Analysis and Scientific Computing, Wiedner Hauptstrasse 8-10, Vienna 1040, Austria}
}

\markboth{}
{}

\IEEEtitleabstractindextext{
\begin{abstract}
The 2D/1D multiscale finite element method (MSFEM) is an efficient way to simulate rotating machines in which each iron sheet is exposed to the same field. It allows the reduction of the three dimensional sheet to a two dimensional cross-section by resolving the dependence along the thickness of the sheet with a polynomial expansion. This work presents an equilibrated error estimator based on flux equilibration and the theorem of Prager and Synge for the T-formulation of the eddy current problem in a 2D/1D MSFEM setting. The estimator is shown to give both a good approximation of the total error and to allow for adaptive mesh refinement by correctly estimating the local error distribution.
\end{abstract}
}

\maketitle

\section{Introduction}
\IEEEPARstart{T}{he} simulation of eddy currents in electrical machines consisting of many steel sheets with the finite element method quickly leads to infeasibly large equation systems. In many machines each sheet in the active zone is exposed to the same field, which allows for a great reduction in computational effort by simulating only a single sheet. However, this reduced problem is still far from trivial. One method to further simplify the problem while maintaining a good approximation of the solution is by spacial decomposition.

The thickness of one sheet is less than a millimeter while the length and width are in the range of meters. A method to treat the two dimensional (2D) cross section and the one dimensional (1D) thickness of the sheet as two coupled problems has been presented in \cite{BottChiamp:02}. It solves the two problems iteratively until convergence is reached. The nature of this coupling has been analyzed in more detail in \cite{PippBelaDlaArkk:10}. 

In \cite{GysSabDul:06} and \cite{HenSteHamGeu:15} different approaches have been presented which isolate the one dimensional problem as a pre-processing step in order to obtain parameters for the two dimensional one.

The 2D/1D multiscale finite element method (MSFEM) presented in \cite{Schoeb:19} uses ideas from the multiscale finite element method to use classic finite element functions for the two dimensional problem while approximating the dependence on the third axis with pre-defined polynomial shape functions, similar to the method presented in \cite{RasDlaFonPippBelaArkkio:11} which is based on trigonometrical shape functions. This enables the solution of the problem within a single iteration while requiring only a mesh for the two dimensional cross section of the sheet. It is also able to include the insulation layers between sheets and correctly treat the edge effect \cite{HollEE:20}.

This paper presents an error estimator for the $\vecT$-formulation of the 2D/1D MSFEM. It is based on flux equilibration and based on the same theory as the error estimator for the $\vecT$-formulation for the MSFEM presented in \cite{Schoeb:21}. In order to fit within the 2D/1D MSFEM framework it has been restructured so both the construction and the evaluation of the estimator require only the two dimensional mesh while being valid in the complete three dimensional domain.

A numerical example shows that the estimator gives a good approximation of the error in both a global and a local sense. This latter property is used to implement adaptive mesh refinement which allows for a high accuracy of the 2D/1D MSFEM solution while requiring significantly less degrees of freedom than uniform mesh refinement.

\section{The $\vecT -\Phi$-Formulation}

We use the $\vecT -\Phi$ formulation for the reference solution of the eddy current problem as described in \cite{Biro:99}. The problem domain $\Omega$ is split into the conducting domain $\Omega_c$, consisting of the steel sheet, and the non-conducting domain $\Omega_0$, consisting of the air regions and the insulation layers. The sheet is assumed to be axis-aligned with the cross-section in the $x-y$ plane and the thickness aligned with the $z$ axis. The total thickness of $\Omega$ is $d=d_{Fe}+d_0$ with the thickness of the sheet $d_{Fe}$ and the thickness of the insulation layer $d_0$.

The magnetic field strength $\vecH$ is written as
\begin{align}
\vecH = \vecT - \nabla\Phi + \vecH_{BS}
\end{align}
with the current vector potential $\vecT\in H(\opcurl ,\Omega_c)$ fulfilling $\opcurl\vecT=\vecJ$, the magnetic scalar potential $\Phi\in H^1(\Omega)$ and a prescribed Biot-Savart field $\vecH_{BS}$. The strong formulation of the eddy current problem in the frequency domain is given as
\begin{align}\label{eq:strong}
\opcurl\rho\opcurl\vecT + i\omega\mu(\vecT - \nabla\Phi + \vecH_{BS}) = \veczero ,
\end{align}
where $\rho=\sigma^{-1}$ is the electric resistivity with the electric conductivity $\sigma$, $\mu$ is the magnetic permeability, $\omega = 2\pi f$ with the frequency $f$ and $i$ is the imaginary unit.

Multiplication with a test function and integration by parts, together with the auxiliary condition $\opdiv\vecB =0$, lead to the weak formulation: Find $\vecT\in H(\opcurl ,\Omega_c)$ and $\Phi\in H^1(\Omega)$ so that
\begin{align}\label{eq:weak3D}
\begin{split}
\int_\Omega \rho\opcurl\vecT\cdot\opcurl\vecv + i\omega\mu(\vecT - \nabla\Phi)\cdot(\vecv - \nabla q) \\
=-\int_\Omega i\omega\mu\vecH_{BS}\cdot(\vecv - \nabla q)
\end{split}
\end{align}
for all $\vecv\in H(\opcurl ,\Omega_c)$ and all $q\in H^1(\Omega)$.

\section{The 2D/1D MSFEM $\vecT$-Formulation}

This paper uses the 2D/1D MSFEM approach for the $\vecT$-formulation which has been described in detail in \cite{Schoeb:19} and \cite{Holl:20}. The three dimensional unknown components $\vecT - \nabla\Phi$ are approximated by
\begin{align}\label{eq:T2D1D}
\vecT_{\text{2D/1D}}= \left( \begin{matrix} 
\phi_0(z)\vecT_{0,x}(x,y) + \phi_2(z)\vecT_{2,x}(x,y) \\ 
\phi_0(z)\vecT_{0,x}(x,y) + \phi_2(z)\vecT_{2,y}(x,y) 
\\ 0
 \end{matrix}\right),
\end{align}
where $\vecT_2\in H(\opcurl_{2D},\Omega_{2D,c})$ and $\vecT_0 = \nabla\Phi_0$ with $\Phi_0\in H^1(\Omega_{2D})$ are defined on the two dimensional projection $\Omega_{2D}$ of $\Omega$. Here and in the following, coordinates $x$, $y$ or $z$ in the index denote the individual components of a vector-valued function. The shape functions $\phi_0$ and $\phi_2$ are predefined piecewise polynomial of order $0$ and $2$, respectively, see also \appendixname~\ref{app:shape}. The two dimensional rotation operator $\opcurl_{2D}$ of a two dimensional vector function $\vecV = (\vecV_x(x,y),\vecV_y(x,y))^T$ is defined as
\begin{align}
\opcurl_{2D}\vecV := \dfrac{\partial}{\partial x}\vecV_y - \dfrac{\partial}{\partial y}\vecV_x.
\end{align}

The discretization of the space $H(\opcurl_{2D})$ is discussed in detail in \cite{Zagl:05}.

The full magnetic field strength is then given by
\begin{align}
\vecH_{\text{2D/1D}} = \vecT_{\text{2D/1D}} + \vecH_{BS} .
\end{align}

For later reference, the (three dimensional) rotation of $\vecH_{\text{2D/1D}}$ is given by
\begin{align}\label{eq:curlT2D1D}
\opcurl\vecH_{\text{2D/1D}}= \left( \begin{matrix} 
- \phi_2'(z)\vecT_{2,y}(x,y) \\ 
\phi_2'(z)\vecT_{2,x}(x,y) \\ 
 \phi_2 \opcurl_{2D}\vecT_2
 \end{matrix}\right) .
\end{align}

To obtain the weak formulation, (\ref{eq:T2D1D}) is used in (\ref{eq:weak3D}) for both the trial function and the test function. Note that $\vecH_{\text{2D/1D}}$ only depends on $z$ via the shape functions $\phi_0$ and $\phi_2$, which are known a-priori. Therefore integration over $z$ can be carried out analytically. This yields the weak 2D/1D MSFEM formulation: Find $\vecT_2\in H(\opcurl_{2D},\Omega_{2D,c})$ and $\Phi_0\in H^1(\Omega_{2D})$ so that
\begin{align}\label{eq:weak2D1D}
\begin{split}
\int_{\Omega_{2D}} & \overline{\rho\phi_2^{'2}}\vecT_2\cdot\vecV_2 + \overline{\rho\phi_2^2}\opcurl_{2D}\vecT_2 \opcurl_{2D}\vecV_2\\
+ & i\omega \left( \overline{\mu\phi_0^{2}}\nabla\Phi_0\cdot\nabla q + \overline{\mu\phi_2^{2}}\vecT_2\cdot\vecV_2 \right) \\
+ & i\omega\overline{\mu\phi_0\phi_2} \left( \nabla\Phi_0\cdot\vecV_2 + \vecT_2\cdot\nabla q \right) \\
=-\int_{\Omega_{2D}} & i\omega\left( \overline{\mu\phi_0^{2} }\vecH_{BS}\cdot \nabla q   + \overline{\mu\phi_0\phi_2}\vecH_{BS}\cdot\vecV_2 \right)  
\end{split}
\end{align}
for all $\vecV_2\in H(\opcurl_{2D},\Omega_{2D,c})$ and $q\in H^1(\Omega_{2D})$ where a bar denotes that the respective function has been integrated with respect to $z$.

\section{Error Estimation}

The proposed error estimator is based on the theorem of Prager and Synge and the theory presented in \cite{BraessSchoeb:08}, which can be adapted to obtain the following identity which is the basis for all further calculations:
\begin{align}\label{eq:prager}
\| \opcurl \vecT - \opcurl \vecT_{\text{2D/1D}} \|_\rho^2 + \| \sigma\gamma - \opcurl \vecT \|_\rho^2 = \| \sigma\gamma - \opcurl \vecT_{\text{2D/1D}} \|_\rho^2 ,
\end{align}
where $\vecT$ is the strong solution of the eddy current problem (\ref{eq:strong}) and $\gamma$ an equilibrated flux fulfilling the condition
\begin{align}\label{eq:gammaeq}
\opcurl\gamma = -i\omega\mu(\vecT_{\text{2D/1D}} + \vecH_{BS}).
\end{align}
The energy norm $\|.\|_\rho$ can be interpreted as a measurement for the eddy current losses, i.e. for the current density $\vecJ$ there holds
\begin{align}
\|\vecJ\|_\rho^2 = \int_{\Omega} \rho\vecJ\cdot\vecJ^\ast = \int_{\Omega} \vecE\cdot\vecJ^\ast,
\end{align}
where the asterisk denotes the complex conjugate.

A variant of (\ref{eq:prager}) for the two dimensional scalar $\vecT$-formulation has been proven in \cite{Schoeb:21} and for the vector-valued magnetostatic case in \cite{BraessSchoeb:08}. The proof of (\ref{eq:prager}) is analogous.

Note that the first term on the left hand side of (\ref{eq:prager}) is the error of the 2D/1D MSFEM solution measured in the norm of the eddy current losses. Assuming a suitable $\gamma$ is known, the right hand side of (\ref{eq:prager}) can be calculated. Given that all terms on the left hand side are guaranteed to be non-negative, the right hand side provides an upper bound for the error.

The main problem is the construction of a suitable $\gamma$ which needs to fulfill (\ref{eq:gammaeq}) on $\Omega$ while at the same time being able to be constructed using only $\Omega_{2D}$. If the error estimator required the full three dimensional domain $\Omega$, it would be much more computationally expensive than the calculation of $\vecT_{2D/1D}$ and nullify the advantages of using a 2D/1D MSFEM. Similarly, the evaluation of the estimator, as defined by the three dimensional integral on the right hand side of (\ref{eq:prager}), needs to be doable using only $\Omega_{2D}$.

This is achieved by using a 2D/1D MSFEM approach for $\gamma$ as well. More specifically, we set
\begin{align}\label{eq:gamma}
\gamma= \left( \begin{matrix} 
\hat\phi_1(z)\gamma_{1,x}(x,y) + \hat\phi_3(z)\gamma_{3,x}(x,y) \\ 
\hat\phi_1(z)\gamma_{1,y}(x,y) + \hat\phi_3(z)\gamma_{3,y}(x,y) \\
\phi_0(z)\gamma_{0,y}(x,y) + \phi_2(z)\gamma_{2,y}(x,y)
 \end{matrix}\right),
\end{align}
with
\begin{align}
\hat\phi_1 := \int \phi_0\, dz  \text{ ~~ and  ~~ } \hat\phi_3 := \int \phi_2\, dz
\end{align}
and the unknowns $\gamma_0,\gamma_2\in H^1(\Omega_{2D})$ and $\gamma_1,\gamma_3\in H(\opcurl_{2D},\Omega_{2D})$ to be determined.

Note that the estimator on the right hand side of (\ref{eq:prager}) consists of $\sigma\gamma$, which is equal to zero in the insulation because of the conductivity, and $\opcurl\vecT_{2D/1D}$, which has only components containing the shape function $\phi_2$, see (\ref{eq:curlT2D1D}), which is also zero in the insulation. Therefore it suffices to consider only the domain of the conducting material for the construction $\gamma$, i.e. $\gamma_0,\gamma_2\in H^1(\Omega_{2D,c})$ and $\gamma_1,\gamma_3\in H(\opcurl_{2D},\Omega_{2D,c})$. A consequence of this is, that $\phi_2' = K\hat\phi_1$ holds with the constant $K = 2\sqrt{6}/d_{Fe}^2$.

The rotation of $\gamma$ is given by
\begin{align}\label{eq:curlgamma}
\opcurl\gamma= \left( \begin{matrix} 
\phi_0\dfrac{\partial}{\partial y}\gamma_0 - \phi_0\gamma_{1,y} + \phi_2\dfrac{\partial}{\partial y}\gamma_2 - \phi_2\gamma_{3,y} \\ 
\phi_0\gamma_{1,x} - \phi_0\dfrac{\partial}{\partial x}\gamma_0 + \phi_2\gamma_{3,x} - \phi_2\dfrac{\partial}{\partial x}\gamma_2 \\
\hat\phi_1\opcurl_{2D}\gamma_1 + \hat\phi_3\opcurl_{2D}\gamma_3
 \end{matrix}\right) .
\end{align}

Writing out the condition (\ref{eq:gammaeq}) using both (\ref{eq:T2D1D}) and (\ref{eq:curlgamma}) and comparing the coefficients with respect to the shape functions yields the equations
\begin{align}
\dfrac{\partial}{\partial y}\gamma_0 - \gamma_{1,y} &= -i\omega\mu (\vecT_{0,x} + \vecH_{BS,x}) , \\
\dfrac{\partial}{\partial y}\gamma_2 - \gamma_{3,y} &= -i\omega\mu\vecT_{2,x} , \\
\gamma_{1,x} - \dfrac{\partial}{\partial x}\gamma_0 &= -i\omega\mu (\vecT_{0,y} + \vecH_{BS,y}) , \\
\gamma_{3,x} - \dfrac{\partial}{\partial x}\gamma_2 &= -i\omega\mu\vecT_{2,y} , \\
\opcurl_{2D}\gamma_1 &= 0 \label{eq:curlgamma1} , \\
\opcurl_{2D}\gamma_3 &= 0 \label{eq:curlgamma3} .
\end{align}

From (\ref{eq:curlgamma1}) and (\ref{eq:curlgamma3}) it follows that $\gamma_1 = \nabla \Phi_1$ and $\gamma_3 = \nabla \Phi_3$ with $\Phi_1,\Phi_3\in H^1(\Omega_{2D,c})$. With this the remaining equations can be rewritten as
\begin{align}
\left( \begin{matrix}  \dfrac{\partial}{\partial y}\gamma_0 - \dfrac{\partial}{\partial y}\Phi_1 \\ \dfrac{\partial}{\partial x}\Phi_1 - \dfrac{\partial}{\partial x}\gamma_0    \end{matrix}\right) &= -i\omega\mu \left( \vecT_{0} + \vecH_{BS} \right) \label{eq:gradgamma1} , \\
\left( \begin{matrix}  \dfrac{\partial}{\partial y}\gamma_2 - \dfrac{\partial}{\partial y}\Phi_3 \\ \dfrac{\partial}{\partial x}\Phi_3 - \dfrac{\partial}{\partial x}\gamma_2    \end{matrix}\right) &= -i\omega\mu \vecT_{2} \label{eq:gradgamma2} .
\end{align}

Note that (\ref{eq:gradgamma1}) and (\ref{eq:gradgamma2}) do not uniquely define all components of $\gamma$. Every solution yields a valid error estimator, but the overestimation (given by the second term on the left hand side of (\ref{eq:prager})) may become arbitrarily large. As can be seen from (\ref{eq:prager}), because the error is independent of $\gamma$, minimizing the overestimation is equivalent to minimizing the estimator. For this purpose additional conditions are imposed. 

Because the estimator is small if $\sigma\gamma$ is a good approximation of $\opcurl\vecT_{\text{2D/1D}}$, a comparison of (\ref{eq:curlT2D1D}) and (\ref{eq:gamma}) suggests that
\begin{align}
\nabla \Phi_1 &\approx K\rho \left( \begin{matrix}  -\vecT_{2,y} \\ \vecT_{2,x}   \end{matrix}\right) , \\
\nabla \Phi_3 &\approx \veczero , \\
\gamma_0 &\approx 0 , \\
\gamma_2 &\approx \rho\opcurl_{2D}\vecT_2
\end{align}
should hold.

In this paper we solve (\ref{eq:gradgamma1}) under the constraint
\begin{align}\label{eq:min1}
\|\sigma\gamma_0\|_\rho^2 + \left\|\sigma\hat\phi_1\nabla\Phi_1 - \phi_2'\left( \begin{matrix}  -\vecT_{2,y} \\ \vecT_{2,x}   \end{matrix}\right)\right\|_\rho^2 \ \to \min
\end{align}
and (\ref{eq:gradgamma2}) under the constraint
\begin{align}\label{eq:min2}
\|\sigma\phi_2\gamma_2 - \phi_2\opcurl_{2D}\vecT_2\|_\rho^2 + \|\sigma\hat\phi_3\nabla\Phi_3\|_\rho^2 \ \to \min .
\end{align}

While this does not yield the optimal minimizer of the estimator because the interdependence of the components is neglected, the numerical example shows that this suffices to achieve an acceptable amount of overestimation. The main advantage of this approach is that instead of one big minimization problem one only has to solve two smaller ones, which is both faster in itself and can even be done in parallel.

The weak formulation for the problem (\ref{eq:gradgamma1}) and (\ref{eq:min1}) reads as: Find $\gamma_0,\Phi_1\in H^1(\Omega_{2D,c})$ and a Lagrange multiplier $\lambda_1\in H_\lambda(\Omega_{2D,c})$ so that
\begin{align}\label{eq:weakest1}
\begin{split}
\int_{\Omega_{2D,c}} & \overline{\sigma\phi_0^2}\gamma_0\chi_0 + \overline{\sigma\hat\phi_1^2}\nabla\Phi_1\cdot\nabla\chi_1 \\
 + & \lambda_1 \cdot \left( \begin{matrix}  \dfrac{\partial}{\partial y}\chi_0 - \dfrac{\partial}{\partial y}\chi_1 \\ \dfrac{\partial}{\partial x}\chi_1 - \dfrac{\partial}{\partial x}\chi_0    \end{matrix}\right)  + \left( \begin{matrix}  \dfrac{\partial}{\partial y}\gamma_0 - \dfrac{\partial}{\partial y}\Phi_1 \\ \dfrac{\partial}{\partial x}\Phi_1 - \dfrac{\partial}{\partial x}\gamma_0    \end{matrix}\right) \cdot \kappa_1 \\
=\int_{\Omega_{2D}} & K\overline{\hat\phi_1^2}  \left( \begin{matrix} -\vecT_{2,y}\\ \vecT_{2,x}  \end{matrix}\right) \cdot \nabla\chi_1 - i\omega\mu \left( \vecT_{0} + \vecH_{BS} \right) \cdot\kappa_1
\end{split}
\end{align}
for all $\chi_0,\chi_1\in H^1(\Omega_{2D,c})$ and $\kappa_1\in H_\lambda(\Omega_{2D,c})$, where, according to the de Rham complex, the Lagrange multiplier space $H_\lambda(\Omega_{2D,c})$ is given as the space $H(\opdiv ,\Omega_{2D,c})$ restricted to divergence-free functions.

Similarly, the weak formulation for the problem (\ref{eq:gradgamma2}) and (\ref{eq:min2}) reads as: Find $\gamma_2,\Phi_3\in H^1(\Omega_{2D,c})$ and a Lagrange multiplier $\lambda_2\in H_\lambda(\Omega_{2D,c})$ so that
\begin{align}\label{eq:weakest2}
\begin{split}
\int_{\Omega_{2D,c}} & \overline{\sigma\phi_2^2}\gamma_2\chi_2 + \overline{\sigma\hat\phi_3^2}\nabla\Phi_3\cdot\nabla\chi_3 \\
+ & \lambda_2 \cdot \left( \begin{matrix}  \dfrac{\partial}{\partial y}\chi_2 - \dfrac{\partial}{\partial y}\chi_3 \\ \dfrac{\partial}{\partial x}\chi_3 - \dfrac{\partial}{\partial x}\chi_2    \end{matrix}\right)  + \left( \begin{matrix}  \dfrac{\partial}{\partial y}\gamma_2 - \dfrac{\partial}{\partial y}\Phi_3 \\ \dfrac{\partial}{\partial x}\Phi_3 - \dfrac{\partial}{\partial x}\gamma_2    \end{matrix}\right) \cdot \kappa_2 \\
=\int_{\Omega_{2D}} & \overline{\phi_2^2}  \opcurl_{2D}\vecT_2 \chi_2 - i\omega\mu\vecT_2 \cdot\kappa_2
\end{split}
\end{align}
for all $\chi_2,\chi_3\in H^1(\Omega_{2D,c})$ and $\kappa_2\in H_\lambda(\Omega_{2D,c})$.

Once all components are calculated, the total estimator can be evaluated on the two dimensional mesh as
\begin{align}\label{eq:esteval}
\begin{split}
\| \sigma\gamma -& \opcurl \vecT_{\text{2D/1D}} \|_\rho^2 = \\
\int_{\Omega_{2D,c}} & \overline{\sigma\hat\phi_1^2} \nabla\Phi_1\cdot\nabla\Phi_1^\ast + \overline{\sigma\hat\phi_3^2} \nabla\Phi_3\cdot\nabla\Phi_3^\ast \\
+ &\overline{\sigma\hat\phi_1\hat\phi_3} \left( \nabla\Phi_1\cdot\nabla\Phi_3^\ast + \nabla\Phi_3\cdot\nabla\Phi_1^\ast \right) \\
- &\overline{\hat\phi_1^2} \left( \nabla\Phi_1\cdot\left( \begin{matrix} -\vecT_{2,y}\\ \vecT_{2,x}  \end{matrix}\right)^\ast + \left( \begin{matrix} -\vecT_{2,y}\\ \vecT_{2,x}  \end{matrix}\right)\cdot\nabla\Phi_1^\ast \right) \\
- &\overline{\hat\phi_1\hat\phi_3} \left( \nabla\Phi_3\cdot\left( \begin{matrix} -\vecT_{2,y}\\ \vecT_{2,x}  \end{matrix}\right)^\ast + \left( \begin{matrix} -\vecT_{2,y}\\ \vecT_{2,x}  \end{matrix}\right)\cdot\nabla\Phi_3^\ast \right) \\
+ &\overline{\rho\hat\phi_1^2}\left( \begin{matrix} -\vecT_{2,y}\\ \vecT_{2,x}  \end{matrix}\right) \cdot \left( \begin{matrix} -\vecT_{2,y}\\ \vecT_{2,x}  \end{matrix}\right)^\ast + \overline{\sigma\phi_0^2}\gamma_0\gamma_0^\ast \\
+ & \overline{\sigma\phi_0\phi_2} (\gamma_0\gamma_2^\ast + \gamma_2\gamma_0^\ast ) + \overline{\sigma\phi_2^2}\gamma_2\gamma_2^\ast \\
- & \overline{\phi_0\phi_2} \left(\gamma_0\opcurl\vecT_2^\ast + \opcurl\vecT_2\gamma_0^\ast \right) \\
- & \overline{\phi_2^2} \left(\gamma_2\opcurl\vecT_2^\ast + \opcurl\vecT_2\gamma_2^\ast \right) \\
+ &\overline{\rho\phi_2^2}\opcurl\vecT_2\opcurl\vecT_2^\ast ,
\end{split}
\end{align}

The integrand can also be used locally to identify the finite elements with the highest contribution to the total error.

\section{Numerical Example}

\begin{figure}[ht]
\centerline{\includegraphics[width=7cm]{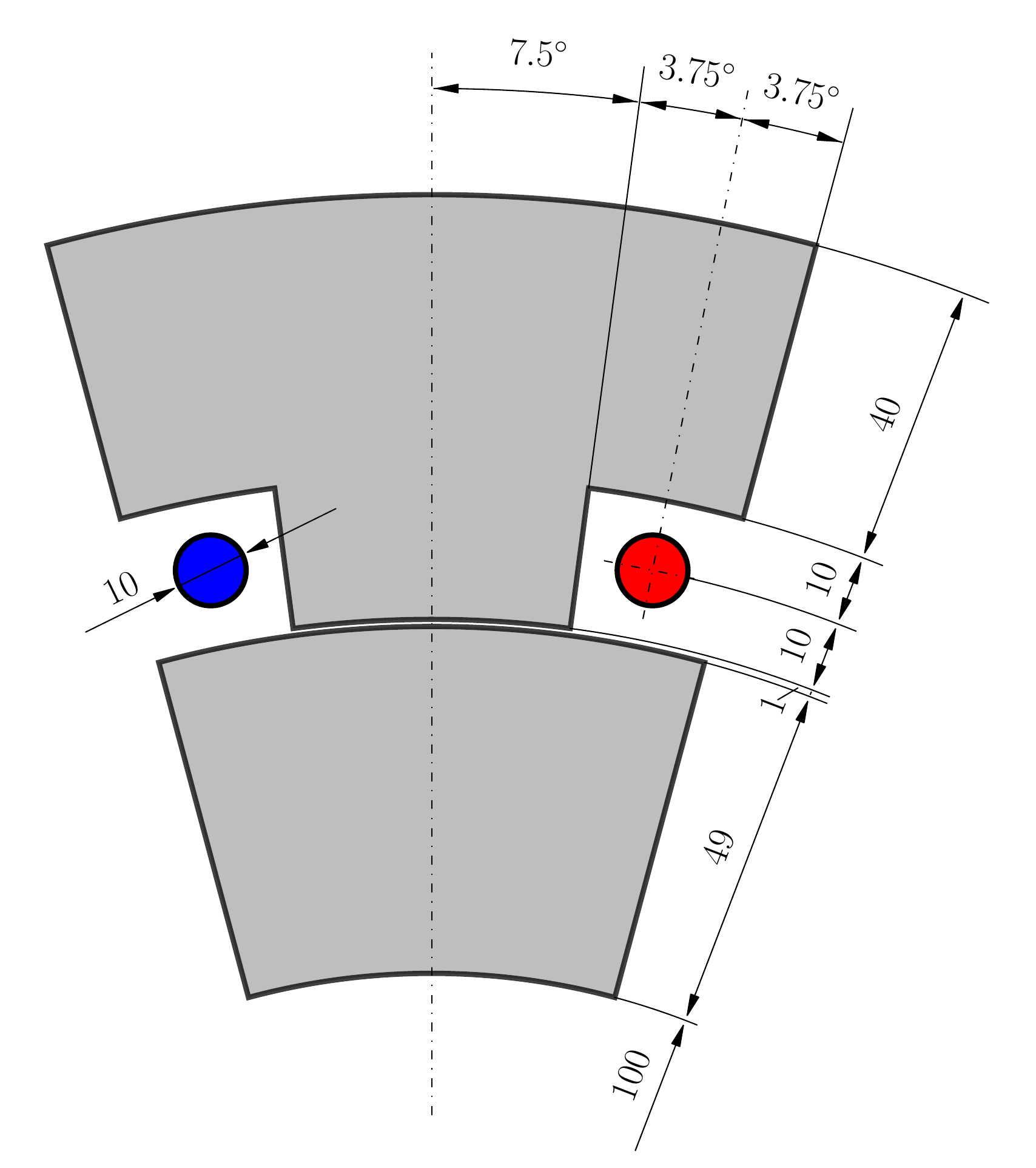}}
\caption{One twelfth of a fictitious machine, consisting of steel sheets (grey) and air domains (left blank). Positive and negative sources are drawn in red and blue, respectively. All measurements are in mm. The thickness of the sheet is $d=0.5$mm with a fill factor of $0.95$.}
\label{fig:domain}
\end{figure}

Consider the machine shown in \figurename~\ref{fig:domain}. Using rotational symmetries, only one twelfth of the entire machine has to be simulated. For the steel sheet a magnetic permeability of $\mu =1000\mu_0$ and an electric conductivity of $\sigma = 2.08$MS is prescribed. The frequency is $50$Hz. The sources are not resolved in the finite element mesh and only included via their Biot-Savart fields. All calculations were done using the open-source software Netgen/NGSolve \cite{NGSolve}.

\begin{figure}[ht]
\centerline{\includegraphics[width=3.5cm]{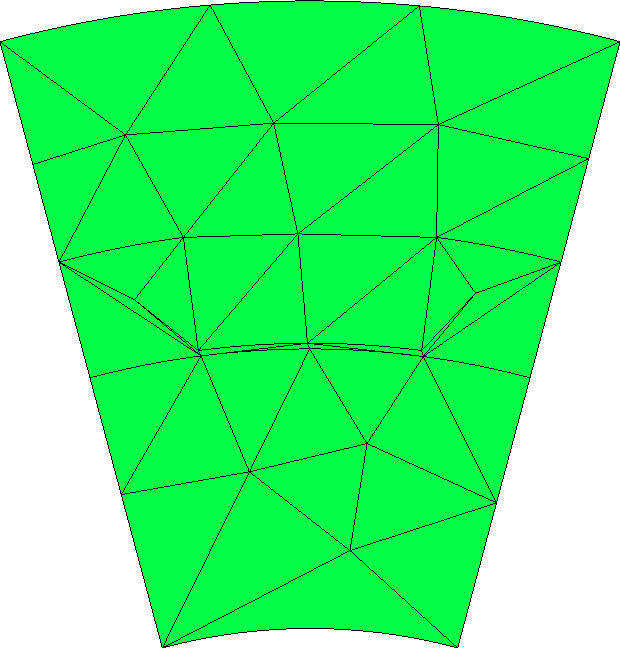}\hspace*{5mm}\includegraphics[width=3.5cm]{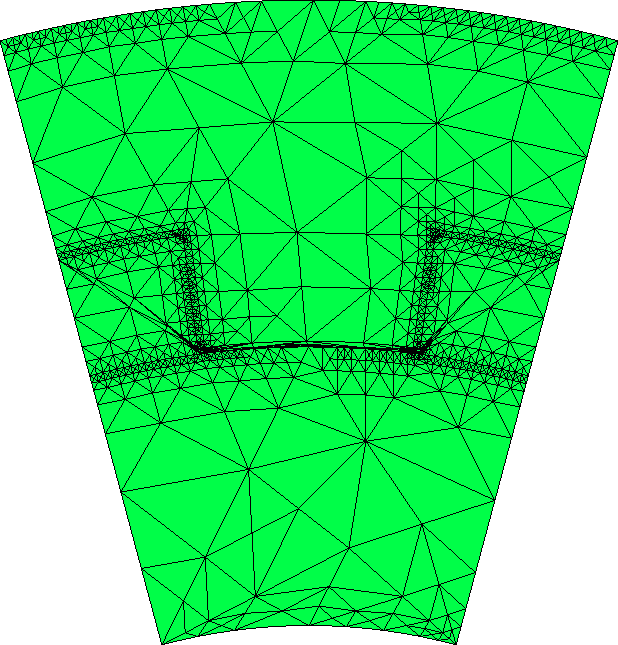}}
\caption{The starting mesh (left) and the adaptively refined mesh after ten iterations (right).}
\label{fig:meshes}
\end{figure}

The calculations start with the coarsest possible mesh for the given geometry, see \figurename~\ref{fig:meshes}. In each iteration, the estimator is evaluated for each individual finite element. Then, all elements where this evaluation yields at least half of the maximum encountered estimator, are refined. In this process some adjacent elements might get refined as well in order to avoid hanging nodes.

\begin{figure}[ht]
\centerline{\includegraphics[width=2.5cm]{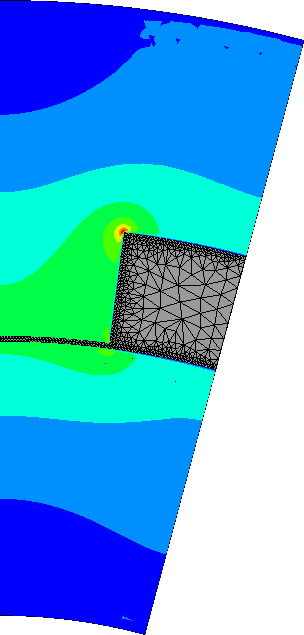}\hspace*{3mm}\includegraphics[width=2.5cm]{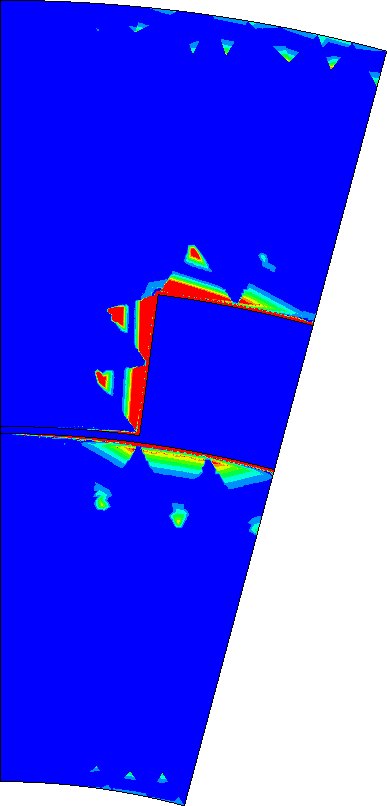}\hspace*{3mm}\includegraphics[width=2.5cm]{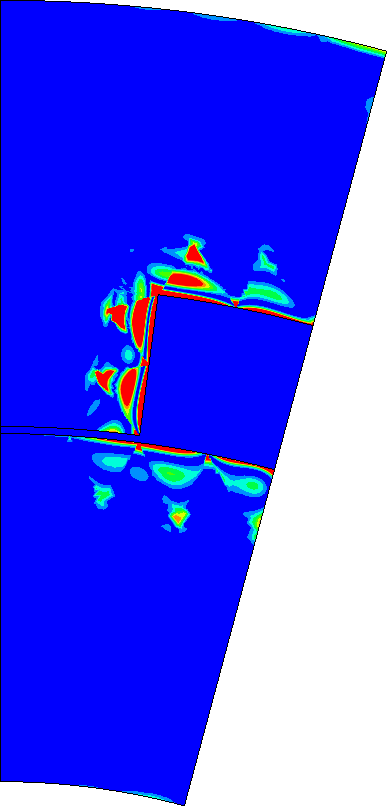}}
\caption{The absolute value of $\vecJ$ drawn for the reference solution (left), the error with respect to the reference solution (middle) and the error estimator (right) after two refinement steps, drawn on one half of the domain.}
\label{fig:Jerrest}
\end{figure}

As can be seen, the refinements are concentrated at the inner edges (where the currents turn around due to the edge effect), at the corners (where the fields peak, see also \figurename~\ref{fig:Jerrest}) and at the inner and outer boundaries (where the boundary conditions need to be resolved correctly). Note also that almost no refinement happens along the vertical symmetry line where the fields are perfectly parallel and easy to resolve.

A qualitative evaluation of the estimator is shown in \figurename~\ref{fig:Jerrest} where both the error (compared to a high order three dimensional reference solution) and the estimator are depicted after two mesh refinements. It can be seen that the estimator correctly identifies the regions where the error is concentrated, further justifying the refinements.

\begin{figure}[ht]
\centerline{\includegraphics[width=7.5cm]{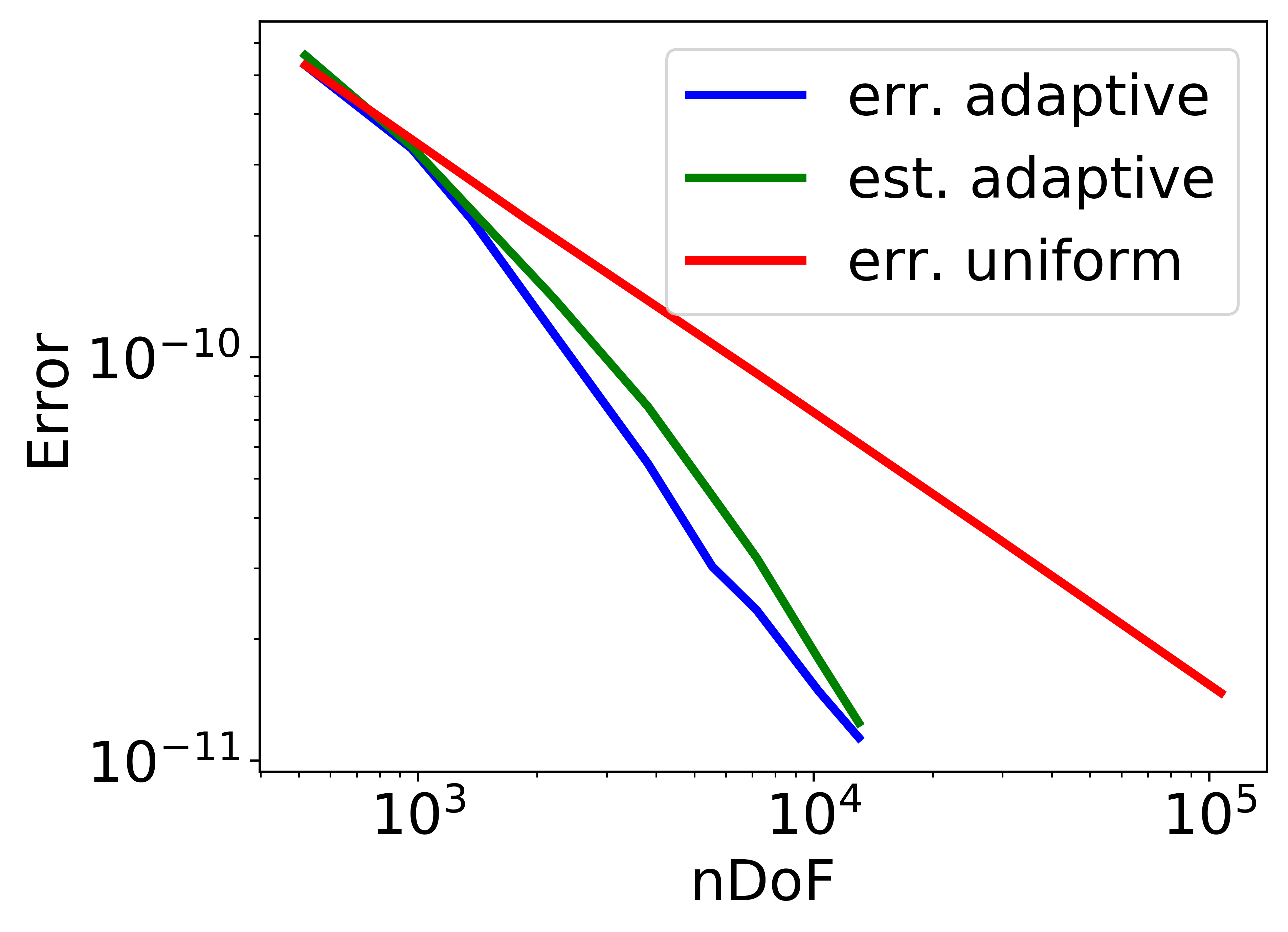}}
\caption{The total error, measured in the norm of the eddy current losses, for both adaptive and uniform mesh refinement as well as the total estimated error for the adaptive approach.}
\label{fig:ndofs}
\end{figure}

Finally, as a quantitative evaluation \figurename~\ref{fig:ndofs} shows the total error of the 2D/1D MSFEM solution compared to the  required degrees of freedom (nDoF) in the finite element problem for both adaptive refinement and uniform refinement. As can be seen, the adaptive refinement leads to a great increase in the rate of convergence. Furthermore, the estimator gives a good approximation of the behavior of the error with only a small overestimation.

\section{Conclusion}
An a-posteriori error estimator has been presented for the 2D/1D MSFEM $\vecT$-formulation of the eddy current problem, based on the theory of flux equlibration. The estimator also utilizes a 2D/1D MSFEM approach in order to prevent its calculation costs to dominate the solution of the problem. Numerical examples show that it gives reliable estimates of the error in both a global and a local sense. This makes it an efficient tool for adaptive mesh refinement to increase the rate of convergence of the 2D/1D MSFEM solution.

\section*{Acknowledgment}
This work was supported by the Austrian Science Fund (FWF) under Project P 31926.

\bibliographystyle{IEEEtran}
\bibliography{BIBLIOGRAPHY}

\begin{thebibliography}{10}
\providecommand{\url}[1]{#1}
\csname url@samestyle\endcsname
\providecommand{\newblock}{\relax}
\providecommand{\bibinfo}[2]{#2}
\providecommand{\BIBentrySTDinterwordspacing}{\spaceskip=0pt\relax}
\providecommand{\BIBentryALTinterwordstretchfactor}{4}
\providecommand{\BIBentryALTinterwordspacing}{\spaceskip=\fontdimen2\font plus
\BIBentryALTinterwordstretchfactor\fontdimen3\font minus
  \fontdimen4\font\relax}
\providecommand{\BIBforeignlanguage}[2]{{%
\expandafter\ifx\csname l@#1\endcsname\relax
\typeout{** WARNING: IEEEtran.bst: No hyphenation pattern has been}%
\typeout{** loaded for the language `#1'. Using the pattern for}%
\typeout{** the default language instead.}%
\else
\language=\csname l@#1\endcsname
\fi
#2}}
\providecommand{\BIBdecl}{\relax}
\BIBdecl

\bibitem{BottChiamp:02}
O.~Bottauscio and M.~Chiampi, ``{Analysis of laminated cores through a directly
  coupled 2-D/1-D electromagnetic field formulation},'' \emph{IEEE Transactions
  on Magnetics}, vol.~38, no.~5, pp. 2358--2360, 2002.

\bibitem{PippBelaDlaArkk:10}
J.~Pippuri, A.~Belahcen, E.~Dlala, and A.~Arkkio, ``{Inclusion of Eddy Currents
  in Laminations in Two-Dimensional Finite Element Analysis},'' \emph{IEEE
  Transactions on Magnetics}, vol.~46, no.~8, pp. 2915--2918, 2010.

\bibitem{GysSabDul:06}
J.~Gyselinck, R.~Sabariego, and P.~Dular, ``{A nonlinear time-domain
  homogenization technique for laminated iron cores in three-dimensional
  finite-element models},'' \emph{IEEE Transactions on Magnetics}, vol.~42,
  no.~4, pp. 763--766, 2006.

\bibitem{HenSteHamGeu:15}
C.~Geuzaine, S.~Steentjes, K.~Hameyer, and F.~Henrotte, ``{P}ragmatic two-step
  homogenisation technique for ferromagnetic laminated cores,'' \emph{IET Meas.
  Sci. Technol.}, vol.~9, no.~2, pp. 152--159, 2015.

\bibitem{Schoeb:19}
M.~Schöbinger, J.~Schöberl, and K.~Hollaus, ``Multiscale fem for the linear
  2-d/1-d problem of eddy currents in thin iron sheets,'' \emph{IEEE
  Transactions on Magnetics}, vol.~55, no.~1, pp. 1--12, 2019.

\bibitem{RasDlaFonPippBelaArkkio:11}
P.~Rasilo \emph{et~al.}, ``{Model of laminated ferromagnetic cores for loss
  prediction in electrical machines},'' \emph{IET Electr. Power Appl.}, vol.~5,
  no.~7, pp. 580--588, 2011.

\bibitem{HollEE:20}
K.~Hollaus and M.~Schöbinger, ``{Air Gap and Edge Effect in the 2-D/1-D Method
  With the Magnetic Vector Potential ${A}$ Using MSFEM},'' \emph{IEEE
  Transactions on Magnetics}, vol.~56, no.~1, pp. 1--5, 2020.

\bibitem{Schoeb:21}
M.~Schöbinger, J.~Schöberl, and K.~Hollaus, ``{An Equilibrated Error
  Estimator for the Multiscale Finite Element Method of a 2-D Eddy Current
  Problem},'' \emph{IEEE Transactions on Magnetics}, vol.~57, no.~6, pp. 1--4,
  2021.

\bibitem{Biro:99}
O.~B{\'{\i}}r{\'o}, ``{Edge element formulations of eddy current problems},''
  \emph{Computer Methods in Applied Mechanics and Engineering}, vol. 169, no.
  3-4, pp. 391--405, 1999.

\bibitem{Holl:20}
K.~Hollaus and M.~Schöbinger, ``{A Mixed Multiscale FEM for the Eddy-Current
  Problem With T, Phi–Phi in Laminated Conducting Media},'' \emph{IEEE
  Transactions on Magnetics}, vol.~56, no.~4, pp. 1--4, 2020.

\bibitem{Zagl:05}
J.~Schoeberl and S.~Zaglmayr, ``{High order Nédélec elements with local
  complete sequence properties},'' \emph{Compel-the International Journal for
  Computation and Mathematics in Electrical and Electronic Engineering -
  COMPEL-INT J COMPUT MATH ELEC}, vol.~24, pp. 374--384, 06 2005.

\bibitem{BraessSchoeb:08}
D.~Braess and J.~Sch{\"o}berl, ``{Equilibrated residual error estimator for
  edge elements},'' \emph{Math. Comp.}, vol.~77, no. 262, pp. 651--672, 2008.

\bibitem{NGSolve}
\BIBentryALTinterwordspacing
J.~Sch{\"o}berl. Netgen/ngsolve. [Online]. Available:
  \url{https://ngsolve.org/}
\BIBentrySTDinterwordspacing

\end{thebibliography}

\appendix

\subsection{Shape Functions}\label{app:shape}

\begin{figure}[ht]
\centerline{\includegraphics[width=7.5cm]{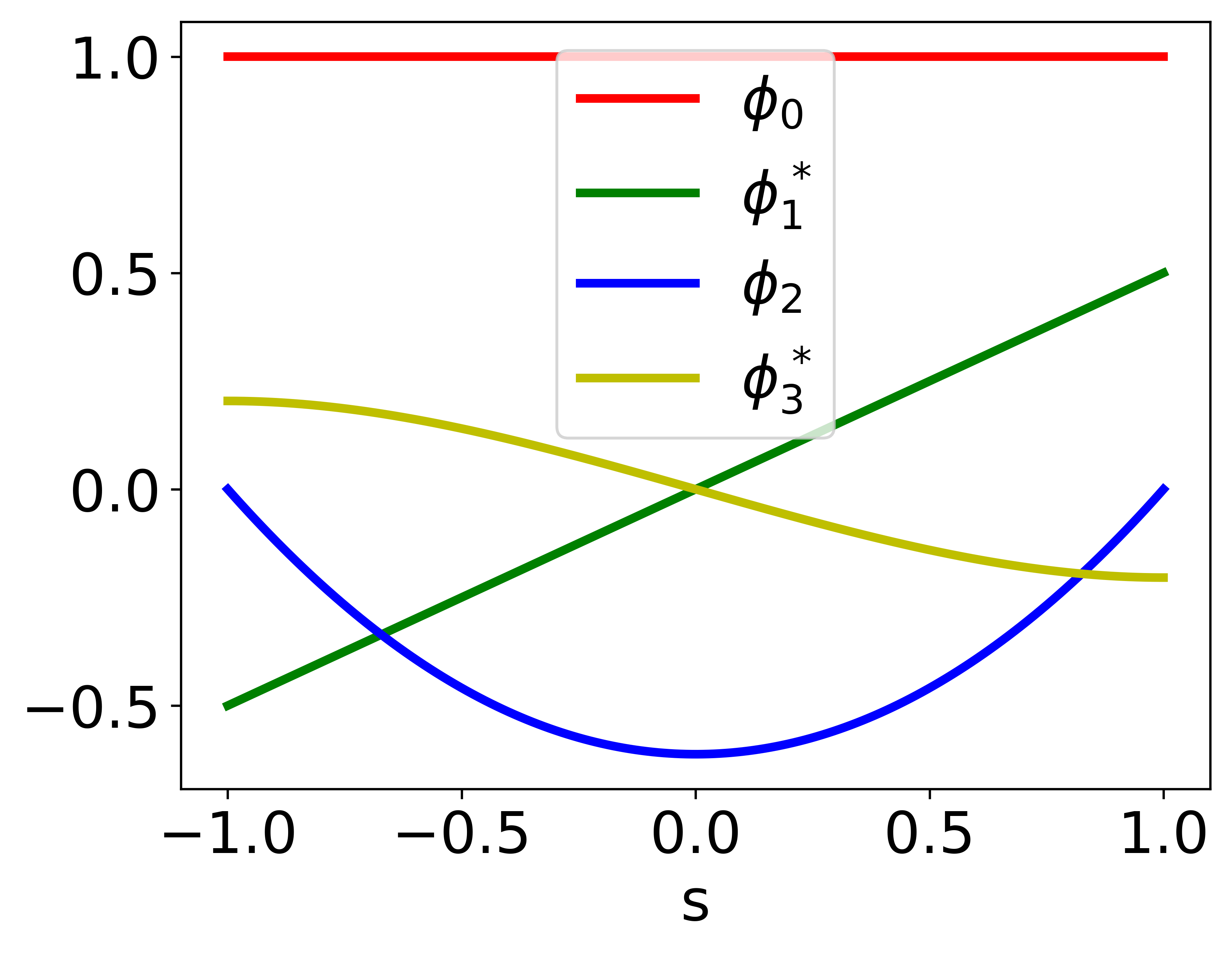}}
\caption{The used shape functions $\phi_0,\hat\phi_1,\phi_2,\hat\phi_3$. Note that $\hat\phi_1$ and $\hat\phi_3$ are scaled to $d_{Fe}=1$ for better visibility.}
\label{fig:phis}
\end{figure}

Assuming that the sheet thickness is aligned with the $z$-axis, the definition of the shape functions uses the auxiliary scaling variable $s:=\frac{2z}{d_{Fe}}$, which transforms the arbitrary interval $\left[ -\frac{d_{Fe}}{2},\frac{d_{Fe}}{2} \right]$ into the normalized interval $[-1,1]$. The shape functions used in this paper are given as
\begin{align}\label{eq:phis}
\phi_0(s) &= 1 , \\
\hat\phi_1(s) &= \frac{d_{Fe}s}{2} ,\\
\phi_2(s) &= \frac{1}{2}\sqrt{\frac{3}{2}}(s^2-1) , \\
\hat\phi_3(s) &= \frac{d_{Fe}\sqrt{6}}{8}s\left( \frac{s^2}{3} -1 \right) ,
\end{align}
see also \figurename~\ref{fig:phis}. In the insulation layer, $\phi_0$ and $\phi_2$ are extended by the constants $1$ and $0$, respectively. The functions $\hat\phi_1$ and $\hat\phi_3$ only appear in the definition for the estimator, which is only defined within the sheet.

In order to write the required integrals concisely, let $\kappa$ be a generic material parameter that is equal to $\kappa_{Fe}$ in the sheet and equal to $\kappa_0$ in the insulation. In the application $\kappa$ takes the place of $\sigma,\rho$ or $\mu$ as needed or it can be omitted by implicitly setting $\kappa_{Fe}=\kappa_0=1$. The integrals required for (\ref{eq:weak2D1D}), (\ref{eq:weakest1}), (\ref{eq:weakest2}) and (\ref{eq:esteval}) are given by
\begin{align}\label{eq:phis}
\overline{\kappa\hat\phi_1^2} &= \frac{d_{Fe}^3\kappa_{Fe}}{12}\\
\overline{\kappa\phi_2^2} &= \frac{d_{Fe}\kappa_{Fe}}{5} \\
\overline{\kappa\phi_2^{'2}} &= \frac{2\kappa_{Fe}}{d_{Fe}} \\
\overline{\kappa\phi_0\phi_2} &= -\frac{\sqrt{6}d_{Fe}\kappa_{Fe}}{6} \\
\overline{\kappa\hat\phi_3^2} &= \frac{17d_{Fe}^3\kappa_{Fe}}{840}\\
\overline{\kappa\hat\phi_1\hat\phi_3} &= -\frac{\sqrt{6}d_{Fe}^3\kappa_{Fe}}{60}
\end{align}

The terms containing $\phi_0^2$ are a special case and one has to differentiate between the integral over the entire domain including the sheet and the insulation and the integral over just the sheet. Therefore
\begin{align}\label{eq:phi0s}
\overline{\kappa\phi_0^2} &= \kappa_{Fe}d_{Fe}+\kappa_0d_0\quad &\text{ in (\ref{eq:weak2D1D})}, \\
\overline{\kappa\phi_0^2} &= \kappa_{Fe}d_{Fe}\qquad &\text{ in (\ref{eq:weakest1}), (\ref{eq:weakest2}) and (\ref{eq:esteval}).}
\end{align}

\end{document}